\documentclass[11pt,letterpaper]{amsart}
\usepackage[parfill]{parskip}
\usepackage{amssymb}
\usepackage{amsmath}
\usepackage{amsthm}
\usepackage{amsfonts}
\usepackage{color}
\usepackage{graphicx, hyperref}
\usepackage{media9}
\usepackage[T1]{fontenc}

\usepackage{textcomp}
\usepackage{palatino}
\usepackage[text={6.5in, 8.5in}, centering]{geometry}

\newtheorem{theorem}{Theorem}

\theoremstyle{definition}

\theoremstyle{remark}
\newtheorem{remark}[theorem]{Remark}

\begin{document}

\title{Balance of the vorticity direction and the vorticity magnitude in 3D fractional Navier-Stokes equations}

\author{Jiayi Wang}
\address{Department of Mathematics, University of Virginia}
\email{jw9yj@virginia.edu}

\date{\today}

\begin{abstract}
Fractional Navier-Stokes equations--featuring a fractional Laplacian--provide a `bridge' between
the Euler equations (zero diffusion) and the Navier-Stokes equations (full diffusion). The problem
of whether an initially smooth flow can spontaneously develop a singularity is a
fundamental problem in mathematical physics, open for the full range of models--from
Euler to Navier-Stokes. The purpose of this work is to present a hybrid, geometric-analytic
regularity criterion for solutions to the 3D fractional Navier-Stokes equations expressed as a
balance--in the average sense--between the vorticity direction and the vorticity magnitude, key geometric
and analytic descriptors of the flow, respectively.
\end{abstract}

\maketitle

\vspace{1in}

\section{Introduction}

The fractional Navier-Stokes equations (NSE) describing the motion of a three-dimensional (3D) Newtonian 
fluid with a reduced diffusion read 

\begin{equation}\label{up}
 \partial_t u + (u \cdot \nabla) u =  - \nu \, ( - \Delta )^\beta u - \nabla p + f
\end{equation}

supplemented with the incompressibility condition $\nabla \cdot u = 0$. 
The diffusion parameter $\beta$ is restricted to the interval $(0, 1)$ (the case 
$\beta = 1$ corresponds to the Navier-Stokes equations).
Here, the vector field $u$ 
is the velocity of the fluid, the scalar field $p$ is the pressure, a positive constant $\nu$ is the viscosity,
and the vector field $f$ is the external force.

\medskip

Henceforth, for simplicity of the exposition, $\nu$ will be set to 1, $f$ taken to 
be a potential force, and the spatial domain will be the whole space.
In this case, the fractional Laplacian, $(-\Delta)^\beta$ is simply a Fourier multiplier 
with the symbol $|\xi|^{2 \beta}$.

\medskip

In order to better understand geometry of the flow mathematically, it is beneficial to
study the vorticity-velocity formulation of the 3D NSE,

\begin{equation}\label{ou}
 \partial_t \omega + (u \cdot \nabla) \omega =  - ( - \Delta)^\beta \omega + (\omega \cdot \nabla) u
\end{equation}

where the vorticity field $\omega$ is given by $\omega = \, \mbox{curl} \, u$. The incompressibility implies 
that
$u$ can be reconstructed from $\omega$ by solving $\Delta u = - \ \mbox{curl} \ \omega$, leading
to the Biot-Savart law, and closing the system. The LHS in (\ref{ou}) is the transport of the vorticity by
the velocity, the first term on the RHS is the fractional diffusion, and the second one is the vortex-stretching
term; there are the three major physical mechanisms in the system.

\medskip

The mathematical theory of the fractional NSE is mostly parallel to the mathematical theory of the Navier-Stokes
equations. It is known that the fractional NSE allow a construction of the Leray-type weak solutions for any $\beta$ in $(0,1)$
(see, e.g., \cite{CoDeLDeR} for a sketch of the proof). In order to be able to apply a standard fixed-point argument and
generate the mild solutions (at least locally-in-time) , one needs the parameter $\beta$ confined to $(\frac{1}{2}, 1)$ (see, e.g., \cite{DL09} where
locally-in-time spatially analytic solutions were constructed). 
Various regularity criteria are also
available. For example, it was shown in \cite{Ch07} that if the (distributional) vorticity $\omega$ of a weak solution
$u$ satisfies 
\begin{equation}\label{x}
 \omega \in L^q(0, T; L^p) \ \ \ \ \ \mbox{with} \ \ \ \ \  \frac{3}{2p}+\frac{\beta}{q} \le \beta
\end{equation}
where $\displaystyle{\frac{3}{\beta} < p \le \infty}$, then the solution is regular on $(0, T]$ (in the endpoint case 
$q=\infty$, in order to avoid smallness in the `standard' argument, the inequality in (\ref{x}) should be strict; a more 
elaborate argument--in the spirit of \cite{ESS03}--allows for the inclusive inequality).

\medskip

In the Navier-Stokes case, a geometric regularity theory--based on coherence of the vorticity field--was 
pioneered in \cite{Co94} where it was shown that the stretching factor in the evolution of the vorticity magnitude
has a singular integral representation featuring a geometric kernel regularized by the local coherence of the
vorticity direction field $\xi$. It was then showed in \cite{CoFe93} that the Lipschitz-coherence of $\xi$ suffices
to rule out singularity formation, establishing the first geometric regularity criterion in the realm
of the `geometric depletion of the nonlinearity'. Among the follow-up works, the Lipschitz
coherence was replaced by the $\frac{1}{2}$-H\"older coherence in \cite{daVeigaBe02}, and
a complete spatiotemporal localization was demonstrated in \cite{Gr09}.

\medskip

It is worth noting that the regularity theory based on the local coherence of the vorticity direction has deep
roots in the computational simulations of turbulent flows which indicate that a dominant morphological
signature of the regions of intense vorticity is the one of the vortex filaments (featuring a high degree of
local coherence of the vorticity direction; see, e.g., \cite{JWSR93, S81, SJO91, VM94}). In addition, long
before the age of high-resolution computational simulations, G.I. Taylor--chiefly based on the experimental
measurements of turbulent flows past a grid--conjectured that stretching of the vortex filaments was the
principal physical mechanism behind the phenomenon of turbulent dissipation (\cite{Tay37}).

\medskip

The article \cite{Ch07} presented a regularity criterion for the fractional NSE in which $\xi$ belongs to a suitable
class of mixed Lebesgue-in temporal variable and Tribel-Lizorkin-in spatial variable spaces while $\omega$
belongs to an interrelated class of the spatiotemporal Lebesgue spaces. A related work \cite{Na19} extended the 
range of some of the functional parameters describing the two classes. The articles were 
in the spirit akin to the results previously obtained in \cite{GrRu04} in the case of the NSE.

\medskip

The regularity criteria in the aforementioned (preceding paragraph) articles quantify a balance between
the vorticity direction and the vorticity magnitude sufficient to rule out a finite time blow-up in the form 
of \emph{two separate}--although interconnected--conditions on $\xi$ and $|\omega|$. A natural question is whether
it is possible to formulate \emph{a single}, hybrid geometric-analytic regularity condition in the form of
a spatiotemporal average 
of the suitably defined degree of coherence of $\xi$ weighted against a power of $|\omega|$
(in the spirit of the NSE article \cite{GrGu10-1}). This would provide a qualitatively rarified measure of the
balance between $\xi$ and $|\omega|$ needed to prevent a possible singularity formation.

\medskip

The goal of this note is to give a positive answer to the above question. More precisely, 
defining a pointwise measure of the coherence of the vorticity direction by

\[
 \rho_\gamma(x,t)=\sup_{y \neq x}
 \frac{|\sin \varphi \bigl(\xi(x+y, t), \xi(x,t)\bigr)|}{|y|^\gamma},
\]

\noindent our main result is summarized as follows.

\begin{theorem}
Let $\omega \in C\bigl( [0, T), L^p\bigr)$ be a (smooth) solution to the fractional NSE (\ref{ou}) 
for some $p > \frac{3}{\beta}$.
Assume that $\omega$
satisfies
\[
 \int_0^T \biggl\{ \int \bigl(\rho_\gamma (x,t) |\omega(x,t)|^a \bigr)^{p_1} \, dx \biggr\}^\frac{2}{p_1} \, dt < \infty
\]
where the parameters $\gamma, p_1$ and $a$ conform to the scaling-invariant condition
\[
 p_1 (\gamma + 2a) - 3 = \beta \, p_1
\]
(in addition to several restrictions to naturally transpire in the proof). Then $T$ is not a blow-up time.
\end{theorem}

\medskip

The next section will first recall some basic concepts in the theory of the geometric depletion of the
nonlinearity and then present the proof of the above theorem.

\bigskip

\section{Balance of $\xi$ and $|\omega|$}

A rigorous study of the geometric depletion of the nonlinearity in
solutions to the 3D NSE is based on a singular integral representation for the stretching
factor in the evolution of the vorticity magnitude $\alpha$ featuring a geometric
kernel $D$ depleted by coherence of the vorticity direction (\cite{Co94}),

\[
 (\partial_t + u \cdot \nabla - \Delta) |\omega|^2 + |\nabla \omega|^2 = \alpha
 |\omega|^2
\]

\noindent where

\[
 \alpha(x, t)=\frac{3}{4\pi} P.V. \int D \left(\hat{y}, \xi(x+y, t),
 \xi(x, t)\right) |\omega(x+y, t)| \frac{1}{|y|^3} \, dy;
\]

\noindent here, $\hat{y}$ is the unit vector in the $y$-direction, $\xi$ is
the vorticity direction and the kernel $D$ is
defined by
\[
 D(e_1, e_2, e_3) = (e_1 \cdot e_3) \left(e_1 \cdot \ (e_2 \times e_3)\right)
\]
for any triple of unit vectors $e_1, e_2$ and $e_3$.

\medskip

Note that

\[
|D \left(\hat{y}, \xi(x+y, t), \xi(x, t)\right)| \le
|\sin\varphi\bigl(\xi(x+y, t), \xi(x, t)\bigr)|
\]

\noindent and--consequently--a coherence condition of the form

\begin{equation}\label{coh}
 |\sin \varphi \bigl(\xi(x+y,t), \xi(x,t)\bigr)| \le c |y|^\delta
\end{equation}

\noindent for some $\delta \in (0, 1)$ will regularize the critical singularity  
of $\displaystyle{\frac{1}{|y|^3}}$ in the integral. This is the essence of the method.

\bigskip

\noindent \emph{Proof of Theorem 1.} \ Multiplying the equations (\ref{ou}) by $\omega |\omega|^{p-2}$
and integrating over the whole space yields

\medskip

\begin{equation}\label{eq}
\frac{1}{p} \frac{d}{dt} \|\omega(t)\|^p_{L^p} + \int \bigl(\sqrt{-\Delta}\bigr)^{2\beta} \omega \cdot \omega |\omega|^{p-2} \, dx
= \int (\omega \cdot \nabla) u \cdot \omega |\omega|^{p-2} \, dx
\end{equation}

\medskip

\noindent (the advection term drops as a result of the divergence-free constraint).

\medskip

Due to the pointwise identity $(\omega \cdot \nabla) u \cdot \omega = \alpha |\omega|^2$, the RHS in
(\ref{eq}) is equal to

\medskip

\[
 \int \alpha |\omega|^p \, dx.
\]

\medskip

A positivity lemma (Lemma 3.3) in \cite{N05} implies a lower bound on the fractional diffusion term,

\medskip

\[
 \int \bigl(\sqrt{-\Delta}\bigr)^{2\beta} \omega \cdot \omega |\omega|^{p-2} \, dx
 \ge \frac{2}{p} \int {\Bigl| \bigl(\sqrt{-\Delta}\bigr)^\beta \bigl(|\omega|^\frac{p}{2}\bigr) \Bigr|}^2 \, dx,
\]

while the Sobolev Embedding Theorem for the fractional derivatives yields

\medskip

\[
 \int {\Bigl| \bigl(\sqrt{-\Delta}\bigr)^\beta \bigl(|\omega|^\frac{p}{2}\bigr) \Bigr|}^2 \, dx 
 \ge c_\beta \biggl( \int |\omega|^\frac{3p}{3-2\beta} \, dx\biggr)^\frac{3-2\beta}{3}
 = c_\beta \, \|\omega(t)\|_{L^\frac{3p}{3-2\beta}}^p.
\]

\medskip

Collecting all of the above, we arrive at the following differential inequality on $(0, T)$,

\medskip

\begin{equation}\label{ineq}
\frac{d}{dt} \|\omega(t)\|^p_{L^p} +  \|\omega(t)\|_{L^\frac{3p}{3-2\beta}}^p
\le c_{p, \beta}  \int \alpha |\omega|^p \, dx.
\end{equation}

\medskip

It will be convenient to estimate the RHS as follows,

\medskip

\begin{equation}\label{split1}
c_{p, \beta} \biggl| \int \alpha(x, t) |\omega(x, t)|^p \, dx \biggr|
\le c_{p, \beta} \, \int \bigl(\rho_\gamma(x, t) |\omega(x, t)|^a\bigr) \biggl( \int \frac{1}{|y|^{3-\gamma}} |\omega(x+y, t)| \, dy \biggr) 
|\omega(x, t)|^b \, dx
\end{equation}

\medskip

\noindent where $a+b=p$.

\medskip

Applying the H\"older inequality with the exponents $\displaystyle{\frac{1}{p_1}+\frac{1}{p_2}+\frac{1}{p_3}=1}$, the quantity
above can be bounded by

\medskip

\begin{equation}\label{split2}
c_{p, \beta} \, \bigl\|\rho_\gamma |\omega|^a\bigr\|_{L^{p_1}} \, \Bigl\| \frac{1}{|y|^{3-\gamma}} \ast |\omega|\Bigr\|_{L^{p_2}}
\, \bigl\||\omega|^b\bigr\|_{L^{p_3}}.
\end{equation}

\medskip

A strong form of the Young inequality for convolutions (estimating one of the functions in the weak Lebesgue space) 
yields the following bound on the middle norm,

\medskip

\[
 \Bigl\| \frac{1}{|y|^{3-\gamma}} \ast |\omega|\Bigr\|_{L^{p_2}}
 \le c \|\omega\|_{L^s}
\]

\noindent where $\displaystyle{\frac{1}{p_2}+1=\frac{3-\gamma}{3}+\frac{1}{s}}$.

\medskip

In addition, the last norm in (\ref{split2}) can be rewritten as

\medskip

\[
 \bigl\||\omega|^b\bigr\|_{L^{p_3}} = \|\omega\|^b_{L^{b  p_3}}.
\]

\medskip

Summarizing, the above estimates on the RHS of (\ref{ineq}) imply

\medskip

\begin{equation}\label{ineqq}
\frac{d}{dt} \|\omega(t)\|^p_{L^p} +  \|\omega(t)\|_{L^\frac{3p}{3-2\beta}}^p
\le c_{p, \beta}  \,   \bigl\|\rho_\gamma |\omega|^a\bigr\|_{L^{p_1}}  \,  \|\omega\|_{L^s} \, \|\omega\|^b_{L^{b  p_3}}.
\end{equation}

\medskip

The next step is to interpolate the last two norms on the RHS between the $L^p$ and the $L^\frac{3p}{3-2\beta}$-norms
of $\omega$,

\medskip

\[
 \|\omega\|_{L^s} \le \|\omega\|_{L^p}^{\alpha'} \, \Bigl\|\omega\Bigr\|_{L^\frac{3p}{3-2\beta}}^{1-\alpha'}
\]

\noindent where $\displaystyle{\frac{1}{s}=\alpha' \frac{1}{p} + (1-\alpha') \frac{3-2\beta}{3p}}$, and

\medskip

\[
 \|\omega\|^b_{L^{b p_3}} \le \|\omega\|_{L^p}^{b \alpha} \, \Bigl\|\omega\Bigr\|_{L^\frac{3p}{3-2\beta}}^{b (1-\alpha)}
\]

\noindent where $\displaystyle{\frac{1}{b p_3}=\alpha \frac{1}{p} + (1-\alpha) \frac{3-2\beta}{3p}}$,

\medskip

\noindent leading to

\medskip

\begin{align*}\label{ineqqq}
\frac{d}{dt} \|\omega(t)\|^p_{L^p} +  \|\omega(t)\|_{L^\frac{3p}{3-2\beta}}^p
& \le c_{p, \beta}  \,   \bigl\|\rho_\gamma |\omega|^a\bigr\|_{L^{p_1}}  \, 
   \|\omega\|_{L^p}^{\alpha' + b \alpha} \, \Bigl\|\omega\Bigr\|_{L^\frac{3p}{3-2\beta}}^{(1-\alpha') + b (1-\alpha)}\\
& \le \widetilde{c}_{p, \beta} \, \bigl\|\rho_\gamma |\omega|^a\bigr\|^2_{L^{p_1}}  \, 
   \|\omega\|_{L^p}^{2\bigl(\alpha' + b \alpha\bigr)} 
   + \frac{1}{2} \Bigl\|\omega\Bigr\|_{L^\frac{3p}{3-2\beta}}^{2 \bigl( (1-\alpha') + b (1-\alpha) \bigr)}.
\end{align*}

\medskip
 
Setting $\displaystyle{ 2\bigl(\alpha' + b \alpha\bigr) = p = 2 \bigl( (1-\alpha') + b (1-\alpha) \bigr)}$ yields the final
form of our differential inequality on $(0, T)$,

\medskip

\begin{equation}\label{ineqqqq}
 \frac{d}{dt} \|\omega(t)\|^p_{L^p}
 \le \widetilde{c}_{p, \beta} \, \bigl\|\rho_\gamma |\omega|^a\bigr\|^2_{L^{p_1}}  \, \|\omega\|^p_{L^p}.
\end{equation}

\medskip

Consequently,

\medskip

\[
  \|\omega(t)\|^p_{L^p} \le  \|\omega_0\|^p_{L^p} \, e^{\Bigl(\widetilde{c}_{p, \beta} \, \int_0^T \bigl\|\rho_\gamma 
  |\omega|^a\bigr\|^2_{L^{p_1}} \, dt\Bigr)}
\]

\medskip

\noindent for all $0 < t < T$.

\medskip

Since $\omega_0 \in L^p$, $\displaystyle{ \sup_{t \in (0, T)} \|\omega(t)\|_{L^p} < \infty}$ provided

\medskip

\[
 \int_0^T \bigl\|\rho_\gamma |\omega|^a\bigr\|^2_{L^{p_1}} \, dt < \infty
\]

\medskip

which--in turn--implies that $T$ is not a blow-up time utilizing the regularity criterion (\ref{x}).
The constraint

\medskip

\[
 p_1 (\gamma + 2a) - 3 = \beta \, p_1
\]

\medskip

is imposed in order to make our hybrid geometric-analytic regularity condition scaling-invariant with respect
to the (unique) intrinsic scaling of the fractional NSE; namely, if a pair $\omega(x, t), \ u(x, t)$ is a solution to
(\ref{ou}), so is the rescaled pair $\omega_\lambda(x, t) = \lambda^{2\beta} \omega\bigl(\lambda x, \lambda^{2\beta} t\bigr),
 \ u_\lambda(x, t) = \lambda^{2\beta-1} u\bigl(\lambda x, \lambda^{2\beta} t\bigr)$, for any $\lambda > 0$.
This type of scaling invariance is useful as it allows one to study the scaling-invariant quantity in view on a
parabolic (corresponding to the fractional diffusion) spatiotemporal cylinder of arbitrary size.
This concludes the proof.

\medskip

\begin{remark}
For a given fractional diffusion parameter $\beta$, there are ten `floating parameters' in the proof, constrained
by eight equations, resulting in a two-parameter family of the hybrid geometric-analytic conditions.
\end{remark}

\begin{remark}
As noted in the introduction, if $\beta > \frac{1}{2}$, the fractional NSE are locally-in-time well posed in $L^p$
for $p$ large enough (cf. \cite{DL09}), and one can replace the assumption
$\omega \in C\bigl( [0, T), L^p\bigr)$ in the theorem simply by an assumption that $\omega_0 \in L^p$ and
designating $T$ to be the first (possible) blow-up time. 
\end{remark}

\bigskip

\bigskip

\bigskip

\centerline{\textbf{Acknowledgments}}

\medskip

The author would like to thank her mentor Professor Zoran Gruji\'c for suggesting the problem and 
his continued guidance, as well as the Department of Mathematics and the Institute of
Mathematical Sciences at the University of Virginia for providing the funding within their
REU-type program in Summer 2019 when the bulk of the project was completed.

\end{document}